\documentclass[reqno]{amsart}

\usepackage[latin1]{inputenc}
\usepackage{amsmath}
\usepackage{amsfonts}
\usepackage{amssymb}
\usepackage{fullpage}
\usepackage{hyperref}
\usepackage{enumitem}
\usepackage{todonotes}
\usepackage{cleveref}
\usepackage{breqn}
\setkeys{breqn}{breakdepth={1}}

\theoremstyle{plain}
\newtheorem{theorem}{Theorem}

\theoremstyle{definition}

\theoremstyle{remark}

\begin{document}
    
    \author{Hieu D. Nguyen}

    \title{A Digital Binomial Theorem}
    \date{12-9-2014}

    \address{Department of Mathematics, Rowan University, Glassboro, NJ 08028.}
    \email{nguyen@rowan.edu}
    
    \subjclass[2010]{Primary 11}
    
    \begin{abstract}
        We present a triangle of connections between the Sierpinski triangle, the sum-of-digits function, and the Binomial Theorem via a one-parameter family of Sierpinski matrices, which encodes a digital version of the Binomial Theorem.
    \end{abstract} 

    \maketitle

\section{Introduction}
It is well known that Sierpinski's triangle can be obtained from Pascal's triangle by evaluating its entries, known as binomial coefficients, mod 2:

\begin{equation}
\begin{array}{c}
1 \\
1 \ \ 1\\
1 \ \ 2 \ \ 1 \\
1 \ \ 3 \ \ 3 \ \ 1 \\
1 \ \ 4 \ \ 6 \ \ 4 \ \ 1 \\
1 \ \ 5 \ 10 \ 10 \ 5 \ \ 1 \\
1 \ \ 6 \ 15 \ 20 \ 15 \ 6 \ \ 1 \\
1 \ \ 7 \ 21 \ 35 \ 35 \ 21 \ 7 \ \ 1 \\
\ldots \\
\textrm{Pascal's triangle}
\end{array}
\begin{array}{c}
\mod 2 \\
 \ \ \ \longrightarrow
\end{array}
\begin{array}{c}
1 \\
1 \ \ 1 \\
1 \ \ 0 \ \ 1 \\
1 \ \ 1 \ \ 1 \ \ 1 \\
1 \ \ 0 \ \ 0 \ \ 0 \ \ 1 \\
1 \ \ 1 \ \ 0 \ \ 0 \ \ 1 \ \ 1 \\
1 \ \ 0 \ \ 1 \ \ 0 \ \ 1 \ \ 0 \ \ 1 \\
1 \ \ 1 \ \ 1 \ \ 1 \ \ 1 \ \ 1 \ \ 1 \ \ 1 \\
\ldots \\
\textrm{Sierpinski's triangle}
\end{array}
\end{equation}
\\
Pascla's triangle is, of course, constructed by inserting the binomial coefficient $\binom{n}{k}$ in the $k$-th position of the $n$-th row, where the first row and first element in each row correspond to $n=0$ and $k=0$, respectively.  Binomial coefficients have a distinguished history and appear in the much-celebrated Binomial Theorem:
\begin{theorem}[Binomial Theorem]
\begin{equation}\label{eq:binomial-theorem}
(x+y)^n = \sum_{k=0}^n \binom{n}{k}x^ky^{n-k}, \ \ \ \ n\in \mathbb{N},
\end{equation}
where $\binom{n}{k}$ are defined in terms of factorials:
\[
\binom{n}{k}=\frac{n!}{k!(n-k)!}.
\]
\end{theorem}
In this article, we demonstrate how the Binomial Theorem in turn arises from a one-parameter generalization of the Sierpinski triangle.  The connection between them is given by the sum-of-digits function, $s(k)$, defined as the sum of the digits in the binary representation of $k$ (see \cite{AS}).  For example, $s(3)=s(1\cdot 2^1+1\cdot 2^0)=2$.  Towards this end, we begin with a well-known matrix formulation of Sierpinski's triangle that demonstrates its fractal nature (see \cite{S}, p.246).  Define a sequence of matrices $S_n$ of size $2^n\times 2^n$ recursively by
\begin{equation}
S_1=\left(
\begin{array}{rr}
1 & 0 \\
1 & 1 \\
\end{array}
\right)
\end{equation}
and
\begin{equation}
 \ \ S_{n+1}=S_1\otimes S_n
 \end{equation}
 for $n>1$. Here, the operation $\otimes$ denotes the Kronecker product of two matrices.  For example, $S_2$ and $S_3$ can be computed as follows:
\begin{align*}
S_2 & =S_1\otimes S_1 = \left(\begin{array}{rr} 1\cdot S_1 & 0\cdot S_1 \\ 1 \cdot S_1 & 1\cdot S_1 \end{array} \right)
= \left(
\begin{array}{rrrr}
1 & 0 & 0 & 0 \\
1 & 1 & 0 & 0 \\
1 & 0 & 1 & 0 \\
1 & 1 & 1 & 1
\end{array}
\right)
\\
\\
S_3 & =S_1\otimes S_2 = \left(\begin{array}{rr} 1\cdot S_2 & 0\cdot S_2 \\ 1\cdot S_2 & 1\cdot S_2 \end{array} \right)
=\left(
\begin{array}{llllllll}
1  \\
1 & 1  \\
1 & 0 & 1  \\
1 & 1 & 1 & 1 \\
1 & 0 & 0 & 0 & 1 \\
1 & 1 & 0 & 0 & 1 & 1 \\
1 & 0 & 1 & 0 & 1 & 0 & 1 \\
1 & 1 & 1 & 1 & 1 & 1 & 1 & 1
\end{array}
\right)
\end{align*}
Thus, in the limit we obtain Sierpinski's matrix $S=\lim_{n\rightarrow \infty}S_n$.  

Less well-known is a one-parameter generalization of Sierpinski's triangle in terms of the sum-of-digits function due to Callan \cite{C}.  If we define
\begin{equation}
S_1(x)=\left(
\begin{array}{rr}
1 & 0 \\
x & 1 \\
\end{array}
\right)
\end{equation}
and
\begin{equation}
 \ \ S_{n+1}(x)=S_1(x) \otimes S_n(x)
 \end{equation}
 for $n>1$, then 
 \begin{equation} \label{eq:generalized-Sierpinski-matrix}
S(x):=\lim_{n\rightarrow \infty} S_n(x) = \left(
\begin{array}{lllllllll}
1  \\
x & 1  \\
x & 0 & 1  \\
x^2 & x & x & 1 \\
x & 0 & 0 & 0 & 1 \\
x^2 & x & 0 & 0 & x & 1 \\
x^2 & 0 & x & 0 & x & 0 & 1 \\
x^3 & x^2 & x^2 & x & x^2 & x & x & 1 \\
\dots & & & & & & & &  \ddots
\end{array}
\right)
\end{equation}
Observe that $S_n(1)=S_n$ and $S(1)=S$.  The matrix $S(x)$ appears  in \cite{C} where Callan defines its entries in terms of the sum-of-digits function $s(k)$.  In particular, if we denote $S(x)=(s_{j,k})$ and assume the indices $j,k$ to be non-negative with $(j,l)=(0,0)$ corresponding to the top left-most entry, then the entries $s_{j,k}$ are defined by
\begin{equation} \label{eq:formula-for-entries-S(x)}
s_{j,k}=\left\{
\begin{array}{cl}
x^{s(j-k)}, & \mathrm{if} \ 0\leq k\leq j \ \textrm{and} \  (k, j-k) \ \textrm{is carry-free}  \\ 
0, & \mathrm{otherwise}
\end{array}
\right. ,
\end{equation}
where the notion of carry-free is defined as follows: call a pair of non-negative integers $(a,b)$ {\em carry-free} if their sum $a+b$ involves no carries when the addition is performed in binary.  For example, the pair $(8,2)$ is carry-free since $8+2=(1\cdot 2^3+0\cdot 2^2+0\cdot 2^1+0\cdot 2^0)+(1\cdot 2^1)=10$ involves no carries in binary.  

To see why (\ref{eq:formula-for-entries-S(x)}) correctly describes (\ref{eq:generalized-Sierpinski-matrix}), we argue by induction.  Clearly, $S_1(x)$ satisfies  (\ref{eq:formula-for-entries-S(x)}).  Next, assume that $S_n(x)$ satisfies  (\ref{eq:formula-for-entries-S(x)}).  It suffices to show that every entry $s_{j,k}$ of $S_{n+1}(x)$ satisfies  (\ref{eq:formula-for-entries-S(x)}).  To prove this, we divide $S_{n+1}(x)$, whose size is $2^{n+1}\times 2^{n+1}$, into four sub-matrices $A,B,C,D$, each of size $2^n\times 2^n$,  based on the recurrence
\[
S_{n+1}(x)= \left( \begin{array}{cc} S_n(x) & 0 \\ xS_n(x) & S_n(x) \end{array}\right)= \left( \begin{array}{cc} A & B \\ C & D\end{array}\right),
\]
where $A=D=S_n(x)$, $B=0$, and $C=xS_n(x)$.  We now consider four cases depending on which sub-matrix the element $s_{j,k}$ belongs to.
\\

\noindent Case 1: $0\leq j,k \leq 2^n-1$.  Then $s_{j,k}$ lies in $A=S_n(x)$ and thus  (\ref{eq:formula-for-entries-S(x)}) clearly holds.
\\

\noindent Case 2: $0\leq j \leq 2^n-1$, $2^n\leq k \leq 2^{n+1}-1$.  Then $s_{j,k}$ lies in $B=0$, which implies $s_{j,k}=0$, and thus   (\ref{eq:formula-for-entries-S(x)})  holds since $k\geq j$.
\\

\noindent Case 3: $2^n \leq j,k \leq 2^{n+1}-1$.  Then $s_{j,k}$ lies in $D=S_n(x)$.  Let
\begin{align*}
j=j_02^0+...+j_{n}2^n, \\
k=k_02^0+...+k_{n}2^n
\end{align*}
denote their binary expansions.  Observe that $j_{n}=k_{n}=1$.  
Define $j'=j-2^n$ and $k'=k-2^n$ where we delete the digit $j_n$ from $j$ (resp. $k_n$ from $k$).  Then it is clear that $(k,j-k)$ being carry-free is equivalent to $(k',j'-k')$ being carry-free.  Moreover, $s(j-k)=s(j'-k')$. We conclude that
\[
s_{j,k}=s_{j',k'}=x^{s(j'-k')}=x^{s(j-k)}
\]
satisfies  (\ref{eq:formula-for-entries-S(x)}).
\\

\noindent Case 4: $2^n \leq j \leq 2^{n+1}-1$, $0\leq k \leq 2^n-1$.  Then $s_{j,k}$ lies in $C=xS_n(x)$.  Define $j'=j-2^n$ and $k'=k$.  Then again $(k,j-k)$ being carry-free is equivalent to $(k',j'-k')$ being carry-free.  Also, $s(j-k)=s(2^n+j'-k')=1+s(j'-k')$.  Hence,
\[
s_{j,k}=xs_{j',k'}=x^{1+s(j'-k')}=x^{s(j-k)}
\]
satisfies  (\ref{eq:formula-for-entries-S(x)}) as well.  This complete the proof.
\\

Callan also proved in the same paper that $S(x)$ generates a one-parameter group, i.e., it satisfies the following additive property under matrix multiplication:
\begin{equation}\label{eq:additive-property}
S(x)S(y)=S(x+y)
\end{equation}
We will see that this property encodes a digital version of the Binomial Theorem.  For example, equating the $(3,0)$-entry of $S(x+y)$, i.e. $s_{3,0}$, with the corresponding entry of $S(x)S(y)$ yields the identity
\begin{equation} \label{eq:sum-of-digit-expansion-n=2}
(x+y)^{s(3)}  = x^{s(3)}y^{s(0)}+x^{s(2)}y^{s(1)}+x^{s(1)}y^{s(2)}+x^{s(0)}y^{s(3)},
\end{equation}
which simplifies to the Binomial Theorem for $n=2$:
\begin{equation}
(x+y)^2  =x^2+2xy+y^2.
\end{equation}
The identities corresponding to the $(5,0)$ and $(7,0)$-entries of $S(x+y)$ are
\begin{equation} \label{eq:digit-binomial-theorem-m=5}
(x+y)^{s(5)}  =  x^{s(5)}y^{s(0)}+x^{s(4)}y^{s(1)}+x^{s(1)}y^{s(4)}+x^{s(5)}y^{s(0)}
\end{equation}
and
\begin{align}  \label{eq:digit-binomial-theorem-m=7}
(x+y)^{s(7)} & = x^{s(7)}y^{s(0)}+x^{s(6)}y^{s(1)}+x^{s(5)}y^{s(2)}+x^{s(4)}y^{s(3)} \\
& \ \ \ \ +x^{s(3)}y^{s(4)}+x^{s(2)}y^{s(5)}+x^{s(1)}y^{s(6)}+x^{s(2)}y^{s(1)}, \notag
\end{align}
respectively.  Observe that (\ref{eq:digit-binomial-theorem-m=5}) is equivalent to (\ref{eq:sum-of-digit-expansion-n=2}) while (\ref{eq:digit-binomial-theorem-m=7}) simplifies to the Binomial Theorem for $n=3$.  

More generally,  property (\ref{eq:additive-property}) can be restated as a digital version of the Binomial Theorem:
\begin{theorem}[Digital Binomial Theorem] \label{th:digital-binomial-theorem}
 Let $m\in\mathbb{N}$.  Then
\begin{equation} \label{eq:digital-binomial-theorem-carry-free}
(x+y)^{s(m)}  = \sum_{\substack{0\leq k\leq m \\ (k,m-k) \ \textrm{carry-free}}}x^{s(k)}y^{s(m-k)}.
\end{equation}
\end{theorem}
We note that (\ref{eq:digital-binomial-theorem-carry-free}) appears implicitly in Callan's proof of (\ref{eq:additive-property}).  The rest of this article is devoted to proving Theorem \ref{th:digital-binomial-theorem} independently of (\ref{eq:additive-property}) and demonstrating that it is equivalent to the Binomial Theorem when $m=2^n-1$.

\section{Proof of the Digital Binomial Theorem}

There are many known proofs of the Binomial Theorem.  The standard combinatorial proof relies on enumerating $n$-element permutations that contain the symbols $x$ and $y$ and then counting those permutations that contain $k$ copies of $x$.  For example, the expansion
\begin{equation} \label{eq:binomial-expansion-n=2-permutation}
(x+y)^2  =xx+xy+yx+yy
\end{equation}
gives all 2-element permutations that contain $x$ and $y$.  Then the number of permutations that contain $k$ copies of $x$ is given by $\binom{2}{k}$.  Thus, (\ref{eq:binomial-expansion-n=2-permutation}) corresponds to (\ref{eq:binomial-theorem}) with $n=2$:
\begin{equation} \label{eq:binomial-theorem-n=2}
(x+y)^2  =\binom{2}{0}x^2+\binom{2}{1}xy+\binom{2}{2}y^2.
\end{equation}

To establish that  (\ref{eq:binomial-theorem-n=2}) is equivalent to (\ref{eq:sum-of-digit-expansion-n=2}), we consider the following {\em digital} binomial expansion: given two sets of digits, $S_0=\{x_0,y_0\}$ and $S_1=\{x_1,y_1\}$, we can represent all ways of constructing a 2-digit  number $z_0z_1$, where $z_0\in S_0$ and $z_1\in S_1$, by the expansion
\begin{equation} \label{eq:digit-expansion-n=2}
(x_0+y_0)(x_1+y_1)  =x_0x_1+x_0y_1+y_0x_1+y_0y_1,
\end{equation}
which we rewrite as 
\begin{equation} \label{eq:digit-expansion-n=2-part2}
(x_0+y_0)(x_1+y_1)  =(x_0^1x_1^1)(y_0^0y_1^0)+(x_0^1x_1^0)(y_0^0y_1^1)+(x_0^0x_1^1)(y_0^1y_1^0)+(x_0^0x_1^0))y_0^1y_1^1).
\end{equation}
If we now assume that $x_0=x_1=x$ and $y_0=y_1=y$, then each term on the right-hand side of (\ref{eq:digit-expansion-n=2-part2}) has the form 
\[
x_0^{d_0}x_1^{d_1} y_0^{1-d_0}y_1^{1-d_1}=x^{s(k)}y^{s(3-k)},
\]
where $k=d_02^0+d_12^1$ and $3-k=(1-d_0)2^0+(1-d_1)2^1$.  It follows that (\ref{eq:digit-expansion-n=2-part2}) reduces to (\ref{eq:sum-of-digit-expansion-n=2}).  On the other hand, (\ref{eq:digit-expansion-n=2}) reduces to (\ref{eq:binomial-expansion-n=2-permutation}).  Thus, we have shown that Theorem \ref{th:digital-binomial-theorem} for $m=3$ is equivalent to the Binomial Theorem for $n=2$.

To extend the proof to integers of the form $m=2^n-1$, we consider $n$ sets of digits, $S_k=\{x_k,y_k\}$, where $k=0,1,\ldots,n-1$.  The expansion
\begin{equation} \label{eq:digital-binomial-expansion}
\prod_{k=0}^{n-1}(x_k+y_k)= \sum_{\substack{z_k\in S_k \\ \forall k=0,1,\ldots,n-1}} z_0\ldots z_{n-1}=\sum_{\substack{d_k \in \{0,1\} \\ \forall k=0,1,\ldots,n-1}} x_0^{d_0}\cdots x_{n-1}^{d_{n-1}} y_0^{1-d_0}\cdots y_{n-1}^{1-d_{n-1}}
\end{equation}
represents all ways of constructing an $n$-digit number $z=z_0z_1\ldots z_{n-1}$ with $z_k\in S_k$ for $k=0,1,\ldots , n-1$.  Then substituting $x_k=x$ and $y_k=y$ for all such $k$ into (\ref{eq:digital-binomial-expansion}) yields
\begin{equation} \label{eq:digital-binomial-expansion-part-2}
(x+y)^n=\sum_{d_0,\ldots,d_{n-1}\in \{0,1\}} x^{d_0+\ldots + d_{n-1}} y^{n-(d_0+\ldots +d_{n-1})},
\end{equation}
or equivalently,
\begin{equation}
(x+y)^{s(2^n-1)}=\sum_{k=0}^{2^n-1} x^{s(k)} y^{s(2^n-1-k)},
\end{equation}
where if we define $k=d_02^0+\ldots + d_{n-1}2^{n-1}$, then $s(k) =d_0+\ldots + d_{n-1}$ and
\[
s(2^n-1-k) =s(2^n-1)-s(k) = n-(d_0+\ldots + d_{n-1}).
\]
Moreover, $k$ ranges from 0 to $2^n-1$ since $d_0,\ldots,d_{n-1}\in \{0,1\}$.  This justifies Theorem 1.  On the other hand, given $k$ between $0$ and $n$, the number of permutations $(d_0,\ldots, d_{n-1})$ containing $k$ 1's is equal to $\binom{n}{k}$.  Thus, (\ref{eq:digital-binomial-expansion-part-2}) reduces to (\ref{eq:binomial-theorem}).  This proves that Theorem 1 is equivalent to the Binomial Theorem.

To complete the proof of Theorem 1 for any non-negative integer $m$, we first expand $m$ in binary: 
\[
m=m_{i_0}2^{i_0}+\ldots +m_{i_{n-1}}2^{i_{n-1}},
\]
where we only record its 1's digits so that $m_{i_k}=1$ for all $k=0,\ldots,n-1$.  Then $s(m)=m_{i_0}+\cdots + m_{i_{n-1}}=n$.  Just as before, we use the expansion (\ref{eq:digital-binomial-expansion}) to derive (\ref{eq:digital-binomial-expansion-part-2}), but this time we rewrite (\ref{eq:digital-binomial-expansion-part-2}) as
\begin{equation}
(x+y)^{s(m)}=\sum_{\substack{0\leq k\leq m \\ (k,m-k) \ \textrm{carry-free}}} x^{s(k)} y^{s(m-k)},
\end{equation}
where we define
\begin{equation}\label{eq:k-expansion}
k=d_02^{i_0}+\ldots + d_{n-1}2^{i_{n-1}}.
\end{equation}
Then $s(k)=d_0+\ldots + d_{n-1}$ and since $m_{i_k}=1$ for all $k=0,\ldots,n-1$, we have
\begin{align*}
m-k & =(m_{i_0}2^{i_0}+\ldots +m_{i_{n-1}}2^{i_{n-1}}) - (d_02^{i_0}+\ldots + d_{n-1}2^{i_{n-1}})) \\
& = (1-d_0)2^{i_0}+\ldots + (1-d_{n-1})2^{i_{n-1}}.
\end{align*}
It follows that
\begin{align*}
s(m-k)  & = (1-d_0)+\ldots + (1-d_{n-1}) \\
& = n-(d_0+\ldots + d_{n-1}).
\end{align*}
Moreover,  it is clear that $0\leq k \leq m$ and $(k,m-k)$ is carry-free.  Conversely, every non-negative integer $k$ with $(k,m-k)$ carry-free must have representation in the form (\ref{eq:k-expansion}); otherwise, the sum $k+(m-k)$ requires a carry in any non-zero digit of $k$ where the corresponding digit of $m$ in the same position is zero.  Thus, Theorem 1 holds for any non-negative integer $m$.

To complete our story we explain why Sierpinski's triangle appears in the reduction of Pascal's triangle's mod 2 by relating binomial coefficients with the sum-of-digits function.  Define the carry function $c(n,k)$ to be the number of carries needed to add $k$ and $n-k$ in binary.  A theorem of Kummer's (see \cite{G}) tells us that the $p$-adic valuation of binomial coefficients is given by the carry function.
\begin{theorem}[Kummer] Let $p$ be a prime integer.  Then the largest power of $p$ that divides $\binom{n}{k}$ equals $c(n,k)$.
\end{theorem}
Kummer's theorem now explains the location of 0's and 1's in Sierpinski's triangle, assuming that its entries are defined by
\begin{equation}\label{eq:Sierpinski-triangle-mod-2}
s_{n,k}=\binom{n}{k} \mod 2.
\end{equation}
Let $p=2$.  If  $(k,n-k)$ is carry free, then $c(n,k)$=0 and therefore the largest power of 2 dividing $\binom{n}{k}$ is $2^0=1$.  In other words, $\binom{n}{k}$ is odd and hence, $s_{n,k}=0$.  On the other hand, if $(k,n-k)$ is \textit{not} carry-free, then $c(n,k)\geq 1$ and so the largest power of 2 dividing $\binom{n}{k}$ is at least 1.  Therefore, $\binom{n}{k}$ is even and hence, $s_{j,k}=0$.  This proves that definition (\ref{eq:Sierpinski-triangle-mod-2}) for Sierpinski's triangle is equivalent to definition (\ref{eq:formula-for-entries-S(x)}) in terms of carry-free pairs with $x=1$. 

Lastly, it is known that the failure of the sum-of-digits function to be additive is characterized by the carry function.  In particular, we have (see \cite{BEJ})
\begin{equation}
s(k)+s(n-k)-s(n)=c(n,k)
\end{equation}
It follows that $(k,n-k)$ is carry-free if and only if $s(k)+s(n-k)=s(n)$.  Thus, it is fitting that the Digital Binomial Theorem can be restated purely in terms of the additivity of the sum-of-digits function:
\begin{equation} \label{eq:digital-binomial-theorem-carry-free}
(x+y)^{s(m)}  = \sum_{\substack{0\leq k\leq m \\ s(k)+s(m-k)=s(m)}}x^{s(k)}y^{s(m-k)}.
\end{equation}


\end{document}